\documentclass{amsart}

\usepackage[cmtip,matrix,arrow]{xy}
\CompileMatrices
\newdir{((}{{\lhook\kern-.1em}}
\newdir{))}{{\rhook\kern-.1em}}
\newdir{ >}{{}*!/-5pt/@{>}}

\usepackage{euler,eucal}

\DeclareMathAlphabet{\mathbf}{T1}{ppl}{bx}{n}
\DeclareMathAlphabet{\mathrm}{T1}{ppl}{m}{n}

\numberwithin{equation}{section}
\renewcommand\dots{\relax\ifmmode\ldots\else$\,\ldots\,$\fi}

\newcommand\note[1]%
{$^\dagger$\marginpar{\footnotesize{$^\dagger${#1}}}}

\divide\count253 by 60 %hours after midnight, rounded down
\divide\count255 by 60 
\multiply\count255 by 60 %minutes after midnight, rounded down
\advance\count254 by -\count255 %minutes after the hour

\def\today{\number\day\space\ifcase\month\or January\or February\or
March\or April\or May\or June\or July\or August\or September\or
October\or November\or December\fi\space\number\year}

\def\hour{\ifnum\count253<10
0\number\count253\else\number\count253\fi}

\def\minute{\ifnum\count254<10
0\number\count254\else\number\count254\fi}

\swapnumbers
\newtheorem{theorem}[equation]{Theorem}
\newtheorem{proposition}[equation]{Proposition}
\newtheorem{lemma}[equation]{Lemma}
\newtheorem{corollary}[equation]{Corollary}

\theoremstyle{definition}

\newtheorem{example}[equation]{Example}

\newcommand\lie{\mathfrak}
\newcommand\g{\lie{g}}
\newcommand\h{\lie{h}} 
 
\newcommand\m{\lie{m}} 
 
\renewcommand\t{\lie{t}}

\newcommand\bb[1]{{\text{\bf#1}}}

\newcommand\R{\bb{R}}

\newcommand\ca{\mathcal}

\newcommand\func[1]{\operatorname{\mathrm{#1}}}
\newcommand\funclim[1]{\operatorname*{\mathrm{#1}}}

\newcommand\Ad{\func{Ad}}

\renewcommand\dim{\func{dim}}

\renewcommand\lim{\funclim{lim}}

\renewcommand\mod{\func{mod}}

\newcommand\pr{\func{pr}}

\newcommand\supp{\func{supp}}

\newcommand\vol{\func{vol}}

\newcommand\Lie{\operatorname{\mathrm{Lie}}} 

\newcommand\group[1]{{\text{\bf#1}}}

\newcommand\U{\group{U}}

\newcommand\abs[1]{\lvert#1\rvert}

\newcommand\quot[1][\kern.3ex]{/\kern-.7ex/_{\kern-.4ex#1}}
\newcommand\bigquot[1][\,\,]{\big/\kern-.85ex\big/_{\!\!#1}}

\newcommand\powl{[\kern-.3ex[}
\newcommand\powr{]\kern-.3ex]}
\newcommand\bigpowl{\bigl[\kern-.6ex\bigl[}
\newcommand\bigpowr{\bigr]\kern-.6ex\bigr]}

\newcommand\inj{\hookrightarrow}
\newcommand\sur{\mathrel{\to\kern-1.8ex\to}}
\newcommand\iso{\mathrel{\hookrightarrow\kern-1.8ex\to}}

\newcommand\longhookrightarrow{\lhook\joinrel\longrightarrow}

\newcommand\longsur{\mathrel{\longrightarrow\kern-1.8ex\to}}
\newcommand\longiso{\mathrel{\longhookrightarrow\kern-1.8ex\to}}

\newcommand\dirac{/\kern-1.2ex\partial} % Dirac operator

\newcommand\antiddots{\mathinner{%
\mkern1mu\raise1pt\vbox{\kern7pt\hbox{.}}
\mkern2mu\raise4pt\hbox{.}\mkern2mu\raise7pt\hbox{.}\mkern1mu}}

\newcommand\inv{^{-1}} 

\renewcommand\subset{\subseteq}

\newcommand\bas{_{\mathrm{bas}}}

\newcommand\prin{_{\mathrm{prin}}}

\begin{document} 

%%%%%%%%%%%%%%%%%%%%%%%%%%%%%%%%%%%%%%%%%%%%%%%%%%%%%%%%%%%%%%%%%%%%%%%%
%%%%%%%%%%%%%%%%%%%%%%%%%%%%%%%%%%%%%%%%%%%%%%%%%%%%%%%%%%%%%%%%%%%%%%%%

\title{A de Rham theorem for symplectic quotients}

\author{Reyer Sjamaar}

\email{sjamaar@math.cornell.edu}

\address{Department of Mathematics, Cornell University, Ithaca, New
York 14853-4201} 

\thanks{The author was partially supported by NSF Grant DMS-0071625}

\date{1 July 2002}

\begin{abstract}
We introduce a de Rham model for stratified spaces arising from
symplectic reduction.  It turns out that the reduced symplectic form
and its powers give rise to well-defined cohomology classes, even on a
singular symplectic quotient.
\end{abstract}

\maketitle

%%%%%%%%%%%%%%%%%%%%%%%%%%%%%%%%%%%%%%%%%%%%%%%%%%%%%%%%%%%%%%%%%%%%%%%%
%%%%%%%%%%%%%%%%%%%%%%%%%%%%%%%%%%%%%%%%%%%%%%%%%%%%%%%%%%%%%%%%%%%%%%%%

%%%%%%%%%%%%%%%%%%%%%%%%%%%%%%%%%%%%%%%%%%%%%%%%%%%%%%%%%%%%%%%%%%%%%%%%
\section{Introduction}
%%%%%%%%%%%%%%%%%%%%%%%%%%%%%%%%%%%%%%%%%%%%%%%%%%%%%%%%%%%%%%%%%%%%%%%%

Let $G$ be a compact Lie group and let $M$ be a smooth $G$-manifold.
Let $\Omega(M)$ be the de Rham complex of differential forms on $M$
and $\Omega\bas(M)$ the subcomplex of basic forms.  It was proved by
Koszul \cite{koszul;certains-groupes} that the cohomology of
$\Omega\bas(M)$ is isomorphic to the cohomology with real coefficients
of the orbit space $M/G$ (which is usually not a manifold, unless $G$
acts freely).

Now suppose that $M$ is equipped with a symplectic form $\omega$ and
that the $G$-action is Hamiltonian with equivariant moment map
$\Phi\colon M\to\g^*$, where $\g=\Lie G$.  The appropriate quotient in
this category is the Marsden-Meyer-Weinstein symplectic quotient
$X=\Phi\inv(0)/G$.  It is usually not a manifold either, unless $G$
acts freely on the fibre $\Phi\inv(0)$, but it always has a natural
stratification into symplectic manifolds.

Much work has been done on the intersection cohomology of symplectic
quotients.  (Cf.\ e.g.\
\cite{kirwan;partial,lerman-tolman;intersection}.)  The purpose of
this note is rather more modest.  We introduce a de Rham model for the
ordinary cohomology of the symplectic quotient $X$, which is a
straightforward adaptation of Koszul's complex of basic forms.  It
relies on a notion of a differential form on $X$ that extends the
concept of a smooth function developed by Arms et al.\
\cite{arms-cushman-gotay;universal}.  Relevant examples are the
reduced symplectic form and its powers, which define cohomology
classes of even degree.  These classes are nonzero if the quotient is
compact.  Thus the symplectic quotient, even when singular, carries a
suitable analogue of a symplectic form and a Liouville volume form.

%%%%%%%%%%%%%%%%%%%%%%%%%%%%%%%%%%%%%%%%%%%%%%%%%%%%%%%%%%%%%%%%%%%%%%%%
\section{Review}\label{section;review}
%%%%%%%%%%%%%%%%%%%%%%%%%%%%%%%%%%%%%%%%%%%%%%%%%%%%%%%%%%%%%%%%%%%%%%%%

Let $(M,\omega)$ be a connected symplectic manifold and let $G$ be a
compact Lie group acting on $M$ in a Hamiltonian fashion with moment
map $\Phi\colon M\to\g^*$, where $\g=\Lie G$.  This means
$d\Phi^\xi=i(\xi_M)\omega$, where $\xi_M$ denotes the vector field on
$M$ induced by $\xi\in\g$ and $\Phi^\xi=\langle\Phi,\xi\rangle$
denotes the component of the moment map along $\xi$.  Also $\Phi$ is
required to be equivariant with respect to the given action on $M$ and
the coadjoint action on $\g^*$.  The \emph{symplectic quotient} of $M$
by $G$ is the space $X=Z/G$, where $Z=\Phi\inv(0)$ is the zero fibre
of the moment map.  It was proved in
\cite{marsden-weinstein;reduction-symplectic-manifolds-symmetry} that
if $G$ acts freely on $Z$, then $Z$ and $X$ are smooth manifolds and
$X$ carries a natural symplectic form.  If $G$ does not act freely on
$Z$, often neither $Z$ nor $X$ are manifolds.  In this case we proceed
as in \cite{sjamaar-lerman;stratified}, the relevant results of which
we recall now.  For any closed subgroup $H$ of $G$ let
$$
M_{(H)}=\{\,m\in M\mid\text{$G_m$ is conjugate to $H$}\,\}
$$
be the stratum of orbit type $(H)$ in the $G$-manifold $M$.  Here
$G_m$ denotes the stabilizer of $m$ with respect to the $G$-action.
Put $Z_{(H)}=Z\cap M_{(H)}$.  Then $Z_{(H)}$ is a smooth $G$-stable
submanifold of $M$.  Let $\{\,Z_a\mid a\in A\,\}$ be the collection of
connected components of all manifolds of the form $Z_{(H)}$, where
$(H)$ ranges over all conjugacy classes of subgroups of $G$.  The
decomposition
\begin{equation}\label{equation;stratification-fibre}
Z=\coprod_{a\in A}Z_a
\end{equation}
is a Whitney stratification of the fibre $Z$.  In particular the index
set $A$ has a partial order defined by $a\le b$ if $Z_a\subset\bar
Z_b$.  There is a unique maximal element in $A$.  The corresponding
stratum, known as the \emph{principal} or \emph{top} stratum $Z\prin$,
is open and dense in $Z$.  Moreover the null-foliation of the
symplectic form $\omega$ restricted to any stratum $Z_a$ is exactly
given by the $G$-orbits.  Hence there exists a unique symplectic form
$\omega_a$ on the quotient manifold $X_a=Z_a/G$ satisfying
$\pi_a^*\omega_a=\iota_a^*\omega$, where $\iota_a\colon Z_a\inj M$ is
the inclusion map and $\pi_a\colon Z_a\sur X_a$ the orbit map.  The
decomposition
\begin{equation}\label{equation;stratification-quotient}
X=\coprod_{a\in A}X_a
\end{equation}
is a locally normally trivial stratification of the quotient $X$ into
the symplectic manifolds $X_a$.  The principal stratum
$X\prin=Z\prin/G$ is open and dense in $X$.

%%%%%%%%%%%%%%%%%%%%%%%%%%%%%%%%%%%%%%%%%%%%%%%%%%%%%%%%%%%%%%%%%%%%%%%%
\section{Forms on a symplectic quotient}\label{section;forms}
%%%%%%%%%%%%%%%%%%%%%%%%%%%%%%%%%%%%%%%%%%%%%%%%%%%%%%%%%%%%%%%%%%%%%%%%

We use the same notation as in the previous section.  We denote the de
Rham complex of a manifold $P$ by $\Omega(P)$.  A \emph{differential
form} on the symplectic quotient $X$ is a differential form $\alpha$
on the top stratum $X\prin$ such that there exists a differential form
$\tilde{\alpha}$ on $M$ satisfying
$\pi\prin^*\alpha=\iota\prin^*\tilde{\alpha}$.  We say that
$\tilde{\alpha}$ \emph{induces} $\alpha$.  An easy averaging argument
shows that we may assume $\tilde{\alpha}$ to be $G$-invariant on $M$.
We denote the collection of differential forms on $X$ by $\Omega(X)$.

If $X=X\prin$, then $X$ and $Z$ are manifolds and the lift of any form
on $X$ to $Z$ can be extended to $M$, so in this case our notion of a
differential form on $X$ reduces to the standard notion.  Observe that
$\Omega(X)$ is a subcomplex of $\Omega(X\prin)$.  Moreover it is
closed under the wedge product.

\begin{example}
The symplectic form $\omega\prin$ on $X\prin$ is induced by the
symplectic form $\omega$ on $M$ and so defines a closed $2$-form on
$X$.
\end{example}

Clearly not every invariant form on $M$ induces a form on $X$.
Indeed, if $\tilde{\alpha}\in\Omega(M)^G$ induces
$\alpha\in\Omega(X)$, then
$\iota\prin^*\tilde{\alpha}=\pi\prin^*\alpha$ is a $G$-horizontal form
on the $G$-manifold $Z\prin$, i.e.\ it is annihilated by all inner
products $i(\xi_M)$ for $\xi\in\g$.  Recall that a form $\beta$ on $M$
is \emph{basic} with respect to the $G$-action if it is $G$-invariant
and $G$-horizontal.  Adapting this notion to our context, let us call
$\beta$ \emph{$\Phi$-basic} if it is $G$-invariant and if
$\iota\prin^*\beta\in\Omega(Z\prin)$ is horizontal.  Let
$\Omega_\Phi(M)$ denote the set of $\Phi$-basic forms.  This is a
subcomplex of $\Omega(M)$ and the kernel of the natural surjection
$\Omega_\Phi(M)\to\Omega(X)$ is the ideal
$$
I_\Phi(M)=\{\,\beta\in\Omega(M)^G\mid\iota\prin^*\beta=0\,\}.
$$
Thus the de Rham complex of $X$ is isomorphic to
\begin{equation}\label{equation;derham}
\Omega(X)\cong\Omega_\Phi(M)/I_\Phi(M),
\end{equation}
a subquotient of the de Rham complex of $M$.  In degree $0$ we have
the smooth functions on $X$ as defined in
\cite{arms-cushman-gotay;universal},
$$
C^\infty(X)\cong C^\infty(M)^G\big/\{\,f\in C^\infty(M)^G\mid
f=0\text{ on }Z\,\}.
$$
If $O$ is a $G$-invariant open neighbourhood of $Z$, then $O$ is a
Hamiltonian $G$-mani\-fold in its own right, so we can define
$\Omega_\Phi(O)$ and $I_\Phi(O)$.  Plainly \eqref{equation;derham}
remains valid if we replace $M$ with $O$.  Thus $\Omega(X)$ depends
only on the $G$-germ of $M$ at $Z$.

It is true, though not completely obvious from the definition, that
every form on $X$ restricts to a form on each stratum of $X$.

\begin{lemma}\label{lemma;restrict}
\begin{enumerate}
\item\label{item;horizontal}
Let $\beta\in\Omega_\Phi(M)$.  Then $\iota_a^*\beta$ is a horizontal
form on $Z_a$ for all $a$.
\item\label{item;null}
Let $\beta\in I_\Phi(M)$.  Then $\iota_a^*\beta=0$ for all $a$.
\item\label{item;restriction}
There is a well-defined restriction map $\Omega(X)\to\Omega(X_a)$ for
each stratum $X_a$.
\end{enumerate}
\end{lemma}

\begin{proof}
Let $\beta\in\Omega_\Phi(M)$ and $z\in Z_a$.  Choose a sequence
$\{z_n\}$ in $Z\prin$ converging to $z$.  By compactness of the
Gra{\ss}mannian we may assume that the sequence of tangent spaces
$T_{z_n}Z\prin$ converges to a subspace $T$ of $T_zM$.  By definition
$i(\xi_M)\beta_{z_n}=0$ on $T_{z_n}Z\prin$ for all $\xi\in\g$, so by
continuity $i(\xi_M)\beta_z=0$ on $T$ for all $\xi$.  By Whitney's
condition (a) we have $T_zZ_a\subset T$.  Hence $i(\xi_M)\beta_z=0$ on
$T_zZ_a$ for all $\xi$.  This proves \eqref{item;horizontal}.

Similarly, if $\beta\in I_\Phi(M)$ then $\beta_{z_n}=0$ on
$T_{z_n}Z\prin$, so by continuity $\beta_z=0$ on $T$ and hence
$\beta_z=0$ on $T_zZ_a$, which proves \eqref{item;null}.

It follows from \eqref{item;horizontal} that if
$\beta\in\Omega_\Phi(M)$ then $\iota_a^*\beta$ descends to a form
$\beta_a$ on $X_a$.  The assignment $\beta\mapsto\beta_a$ defines a
homomorphism $\Omega_\Phi(M)\to\Omega(X_a)$ for each $a$.  It follows
from \eqref{item;null} that this map is $0$ on the ideal $I_\Phi(M)$.
Using the isomorphism \eqref{equation;derham} we obtain the desired
restriction map $\Omega(X)\to\Omega(X_a)$.
\end{proof}

%%%%%%%%%%%%%%%%%%%%%%%%%%%%%%%%%%%%%%%%%%%%%%%%%%%%%%%%%%%%%%%%%%%%%%%%
\section{Symplectic induction}\label{section;induction}
%%%%%%%%%%%%%%%%%%%%%%%%%%%%%%%%%%%%%%%%%%%%%%%%%%%%%%%%%%%%%%%%%%%%%%%%

A shortcoming of the de Rham complex $\Omega(X)$ is that it appears to
depend on the way in which $X$ is written as a quotient.  But in
certain interesting situations this defect turns out to be illusory.
For instance, let $H$ be a closed subgroup of $G$ and let
$(N,\omega_N)$ be a Hamiltonian $H$-manifold with equivariant moment
map $\Psi\colon N\to\h^*$.  Consider the Hamiltonian $G\times H$-space
$$
P=T^*G\times N,
$$
where the action of $G$ on $P$ is given by left multiplication on
$T^*G$ and the action of $H$ by right multiplication on $T^*G$ and the
given action on $N$.  Let $M$ be the symplectic quotient of $P$ with
respect to the $H$-action.  This is called the $G$-space
\emph{induced} by the $H$-space $N$.  Since $H$ acts freely on $T^*G$,
$M$ is a smooth manifold and from $P$ it inherits a symplectic form
$\omega$ and a Hamiltonian $G$-action with moment map $\Phi$.  Let $Y$
be the symplectic quotient of $N$ by the $H$-action and $X$ the
symplectic quotient of $M$ by the $G$-action.  The principle of
reduction in stages implies that $X$ and $Y$ are isomorphic in the
sense that there is a stratification-preserving homeomorphism $Y\to X$
that restricts to a symplectomorphism on each stratum.  We can
represent the situation symbolically by a commutative diagram
$$
\xymatrix{P\ar@{.>}[r]^-{G}\ar@{.>}[d]_{H}
&N\ar@{.>}[d]^{H}\\M\ar@{.>}[r]_-{G}&**[r]Y\cong X,}
$$
where the dotted arrows indicate symplectic reduction with respect to
the relevant group.  We assert that the de Rham complexes of $X$ and
$Y$ are likewise isomorphic.

To prove this we need to recall from
\cite[\S~2]{sjamaar-lerman;stratified} the definition of the
isomorphism $Y\to X$.  Choose an $H$-invariant subspace $\m$ of $\g$
complementary to the subalgebra $\h$.  Then we have $H$-invariant
decompositions $\g=\h\oplus\m$ and $\g^*=\h^*\oplus\m^*$.  Define a
map
$$
G\times\m^*\times N\to P\cong G\times\g^*\times N
$$
by sending $(g,\alpha,p)$ to $(g,\alpha-\Psi(p),p)$.  This is an
$H$-equivariant diffeomorphism from $G\times\m^*\times N$ onto the
zero fibre of the $H$-moment map on $P$.  Taking quotients by $H$
we obtain a $G$-equivariant diffeomorphism
$$
M\cong(G\times\m^*\times N)/H
$$
from $M$ to the homogeneous vector bundle over $G/H$ with fibre
$\m^*\times N$.  Let us identify $M$ with this bundle and write a
typical point in it as $[g,\alpha,p]$ with $g\in G$, $\alpha\in\m^*$
and $p\in N$.  The $G$-action on $M$ is given by
$k[g,\alpha,p]=[kg,\alpha,p]$ for $k\in G$ and the moment map by
\begin{equation}\label{equation;moment}
\Phi([g,\alpha,p])=\Ad^*(g)(\alpha+\Psi(p)).
\end{equation}
Let $f\colon N\to M$ be the embedding defined by $f(p)=[1,0,p]$.  Then
$f$ is $H$-equivariant and \eqref{equation;moment} shows that
$\Phi\circ f=\pr^*\circ\Psi$, where $\pr^*\colon\h^*\to\g^*$ is the
transpose of the projection map $\g\to\h$.  Hence $f$ maps the zero
fibre $Z_N=\Psi\inv(0)$ into $Z$ and descends to a map $Y\to X$, which
is the required isomorphism.  In particular $f$ maps the principal
stratum $(Z_N)\prin$ into the principal stratum of $Z$.  In fact $Z$
and $Z\prin$ are homogeneous bundles over $G/H$,
$$
Z=(G\times Z_N)/H\qquad\text{and}\qquad Z\prin=(G\times(Z_N)\prin)/H.
$$
This implies that the restriction map $f^*\colon\Omega(M)\to\Omega(N)$
sends $\Omega_\Phi(M)$ to $\Omega_\Psi(N)$ and $I_\Phi(M)$ to
$I_\Psi(N)$.  Therefore, because of the isomorphism
\eqref{equation;derham}, it descends to a map
$r\colon\Omega(X)\to\Omega(Y)$.

\begin{proposition}\label{proposition;induction}
The map $r\colon\Omega(X)\to\Omega(Y)$ is an isomorphism.
\end{proposition}

\begin{proof}
This relies on some material developed in Appendix
\ref{section;homogeneous}.  Let $\iota\prin\colon Z\prin\to M$ be the
inclusion map.  This is a bundle map of fibre bundles over the base
$G/H$.  Its restriction to a fibre is the inclusion map
$(\iota_N)\prin\colon(Z_N)\prin\to N$.  Let
\begin{align*}
e_M&\colon\Omega(N)^H\to\Omega(M)^G\\
e_Z&\colon\Omega((Z_N)\prin)^H\to\Omega(Z\prin)^G
\end{align*}
be the extension homomorphisms for the homogeneous bundles $M$ and
$Z\prin$ as defined in Appendix \ref{section;homogeneous}.  Then
\begin{equation}\label{equation;extend}
e_Z\circ(\iota_N)\prin^*=\iota\prin^*\circ e_M
\end{equation}
by Lemma \ref{lemma;functor}.  

Now we can show that $r$ is surjective.  In fact we must show that
$f^*\Omega_\Phi(M)=\Omega_\Psi(N)$.  Let $\gamma\in\Omega_\Psi(N)$.
Then by definition $(\iota_N)\prin^*\gamma$ is $H$-basic, so
$e_Z((\iota_N)\prin^*\gamma)$ is $G$-basic by Lemma
\ref{lemma;extension}\eqref{item;basic}.  From \eqref{equation;extend}
we get that $\iota\prin^*e_M(\gamma)$ is $G$-basic, i.e.\
$e_M(\gamma)\in\Omega_\Phi(M)$.  Using Lemma
\ref{lemma;extension}\eqref{item;right} we find that $\gamma=f^*\beta$
with $\beta=e_M(\gamma)\in\Omega_\Phi(M)$.  Hence
$f^*\Omega_\Phi(M)=\Omega_\Psi(N)$.

Next we prove that $r$ is injective.  Suppose that
$\beta\in\Omega_\Phi(M)$ satisfies $f^*\beta\in I_\Psi(N)$.  We need
to show that $\beta\in I_\Phi(M)$.  The assumptions on $\beta$ mean
that $\iota\prin^*\beta$ is $G$-basic and that
$(\iota_N)\prin^*f^*\beta=0$.  Using Lemma
\ref{lemma;extension}\eqref{item;left} we get
$$
\iota\prin^*\beta=e_Z(f^*\iota\prin^*\beta)
=e_Z((\iota_N)\prin^*f^*\beta)=e_Z(0)=0,
$$
i.e.\  $\beta\in I_\Phi(M)$.
\end{proof}

%%%%%%%%%%%%%%%%%%%%%%%%%%%%%%%%%%%%%%%%%%%%%%%%%%%%%%%%%%%%%%%%%%%%%%%%
\section{The de Rham sheaf}\label{section;sheaf}
%%%%%%%%%%%%%%%%%%%%%%%%%%%%%%%%%%%%%%%%%%%%%%%%%%%%%%%%%%%%%%%%%%%%%%%%

To prove a de Rham theorem we need to sheafify the de Rham complex.
Let $U$ be an open subset of the symplectic quotient $X$.  The
stratification of $X$ induces one on $U$, so we can talk about the
principal stratum of $U$ etc.  A \emph{differential form} on $U$ is a
differential form $\alpha$ on $U\prin$ such that for all $x\in U$
there exist $\alpha'\in\Omega(X)$ and an open neighbourhood $U'$ of
$x$ in $U$ such that $\alpha=\alpha'$ on $U'\prin$.  The set of
differential forms on $U$ is denoted by $\Omega(U)$.  It is easy to
check that the presheaf of differential graded algebras $\Omega\colon
U\mapsto\Omega(U)$ is a sheaf.  Its space of global sections is the
previously defined de Rham complex $\Omega(X)$.

\begin{lemma}
$\Omega$ is an acyclic sheaf, i.e.\ $H^i(X,\Omega^j)=0$ for all
$i\ge1$ and $j\ge0$.
\end{lemma}

\begin{proof}
The space $X$ possesses smooth partitions of unity subordinate to
arbitrary open covers $\ca U$.  Indeed, for each $U\in\ca U$ choose a
$G$-invariant open $\tilde{U}$ in $M$ such that $U=(\tilde{U}\cap
Z)/G$ and let $O$ be the union of the $\tilde{U}$'s.  Choose a smooth
$G$-invariant partition of unity on the $G$-manifold $O$ subordinate
to the cover defined by the $\tilde{U}$'s; this induces a smooth
partition of unity on $X$ subordinate to $\ca U$.  Thus the sheaf of
smooth functions $\Omega^0$ is fine in the sense of
\cite[\S~3.7]{godement;faisceaux}.  A standard result in sheaf theory
(see e.g.\ \cite[Th\'eor\`eme 4.4.3]{godement;faisceaux}) now implies
that $\Omega^0$ is acyclic.  Since $\Omega$ is a module over
$\Omega^0$, it is fine, and therefore acyclic, as well.
\end{proof}

There is an alternative characterization of forms on open subsets of
$X$.  The proof is an easy exercise involving partitions of unity.

\begin{lemma}
Let $U$ be an open subset of $X$ and let $\alpha\in\Omega(U\prin)$.
Then $\alpha\in\Omega(U)$ if and only if there exist a $G$-invariant
open subset $\tilde{U}$ of $M$ and a form
$\tilde{\alpha}\in\Omega(\tilde{U})$ such that $U=(\tilde{U}\cap Z)/G$
and $\iota\prin^*\tilde{\alpha}=\pi\prin^*\alpha$.
\end{lemma}

Now let $\underline{\R}$ be the sheaf of locally constant real-valued
functions on $X$ and consider the sequence
\begin{equation}\label{equation;sequence}
0\to\underline{\R}\xrightarrow{i}\Omega^0\xrightarrow{d}
\Omega^1\xrightarrow{d}\cdots,
\end{equation}
where $i\colon\underline{\R}\to\Omega^0$ is the natural inclusion.
The following assertion is proved in the next section.

\begin{lemma}\label{lemma;poincare}
The sequence \eqref{equation;sequence} is exact.
\end{lemma}

Thus the de Rham complex is an acyclic resolution of the constant
sheaf, which by standard sheaf theory (see e.g.\ \cite[Th\'eor\`emes
4.7.1, 6.2.1]{godement;faisceaux}) implies the following de Rham
theorem.

\begin{theorem}\label{theorem;derham}
The de Rham cohomology ring $H(\Omega(X))$ is naturally isomorphic to
the (\v Cech or singular) cohomology ring of $X$ with real
coefficients $H(X,\R)$.
\end{theorem}

%%%%%%%%%%%%%%%%%%%%%%%%%%%%%%%%%%%%%%%%%%%%%%%%%%%%%%%%%%%%%%%%%%%%%%%%
\section{The Poincar\'e lemma}\label{section;poincare}
%%%%%%%%%%%%%%%%%%%%%%%%%%%%%%%%%%%%%%%%%%%%%%%%%%%%%%%%%%%%%%%%%%%%%%%%

In this section we prove the following (marginally stronger) version
of Lemma \ref{lemma;poincare}: every $x\in X$ has a basis of open
neighbourhoods $U$ such that the sequence
\begin{equation}\label{equation;poincare}
0\to\R\xrightarrow{i}\Omega^0(U)\xrightarrow{d}\Omega^1(U)
\xrightarrow{d}\cdots
\end{equation}
is exact.  The proof is a variation on a familiar homotopy argument in
de Rham theory, which requires a brief look into the functorial
properties of $\Omega(X)$.

Let $(M',\omega',\Phi')$ be a second Hamiltonian $G$-manifold with
zero fibre $Z'=(\Phi')\inv(0)$ and symplectic quotient $X'=Z'/G$.
Then we have stratifications $Z'=\coprod_{a\in A'}Z'_a$ and
$X'=\coprod_{a\in A'}X'_a$ analogous to those for $Z$ and $X$.  Let us
call a map $f\colon M\to M'$ \emph{allowable} if
\begin{enumerate}
\item\label{item;smooth}
$f$ is smooth and $G$-equivariant;
\item\label{item;Z}
$f(Z)\subset Z'$;
\item\label{item;tangent}
$df(T_zZ\prin)\subset T_{f(z)}Z'_{a(z)}$ for all $z\in Z\prin$, where
$Z'_{a(z)}\subset Z'$ is the stratum of $f(z)$.
\end{enumerate}
For instance, if $f$ is smooth and equivariant and maps $Z\prin$ into
a single stratum of $Z'$, then $f$ is allowable.

\begin{example}\label{example;dilation}
Let $(V,\omega)$ be a symplectic vector space on which $G$ acts
linearly and symplectically.  A moment map is given by
$\Phi_V^\xi(v)=\frac1{2}\omega(\xi v,v)$, where $\xi\in\g$ acts on $V$
via the infinitesimal representation $\g\to\lie{sp}(V)$.  Let $t\in\R$
and let $f\colon V\to V$ be the dilation $f(v)=tv$.  Clearly $f$
preserves $Z$.  Furthermore, if $\t\ne0$, then $f(v)$ has the same
stabilizer as $v$, so $f$ maps $Z\prin$ to itself.  If $t=0$, then $f$
maps $Z\prin$ to $0$.  In either case $f$ maps $Z\prin$ into a single
stratum of $Z$ and it is obviously smooth and equivariant, so it is
allowable.  Similarly, if $\abs{t}\le1$ and $B$ is a $G$-invariant
open ball about the origin, the restriction of $f$ is an allowable map
from $B$ to itself.
\end{example}

The following result is easy to deduce from Lemma
\ref{lemma;restrict}.

\begin{lemma}\label{lemma;allowable}
Let $f\colon M\to M'$ allowable.  Then the pullback map
$f^*\colon\Omega(M')\to\Omega(M)$ sends $\Omega_{\Phi'}(M')$ to
$\Omega_\Phi(M)$ and $I_{\Phi'}(M')$ to $I_\Phi(M)$, and therefore
induces a homomorphism $f^*\colon\Omega(X')\to\Omega(X)$.
\end{lemma}

Homotopies induce chain homotopies on the de Rham complex in a
standard way.  Let $F\colon M\times[\,0,1]\to M'$ be a smooth homotopy
and put $F_t=F|_{M\times\{t\}}$.  Let $t$ be the coordinate on
$[\,0,1]$ and for $\gamma\in\Omega(M')$ put
$\kappa_F\gamma=\int_0^1i(\partial/\partial t)F^*\gamma\,dt$.  Then
$\kappa_F$ lowers the degree by $1$ and
$$
F_1^*-F_0^*=\kappa_Fd+d\kappa_F.
$$
Assume that $F$ is equivariant with respect to the given $G$-actions
on $M$ and $M'$ and the trivial action on $[\,0,1]$.  It is
straightforward to check that
\begin{align}
\kappa_F\circ g^*&=g^*\circ\kappa_F&&\text{for all $g\in
G$,}\label{equation;intertwine}\\
\kappa_F\circ i(\xi_{M'})&=-i(\xi_M)\circ\kappa_F&&\text{for all
$\xi\in\g$.}\label{equation;inner}
\end{align}
Call the homotopy $F$ \emph{allowable} if
\begin{enumerate}
\item\label{item;smooth-homotopy}
$F$ is smooth and $G$-equivariant;
\item\label{item;t}
$F_t\colon M\to M'$ is allowable for almost all $t\in[\,0,1]$;
\item\label{item;tangent-homotopy}
$dF_{(z,t)}(\partial/\partial t)\in T_{F(z,t)}Z'_{a(z,t)}$ for almost
all $t\in[\,0,1]$ and for all $z\in Z\prin$, where $Z'_{a(z,t)}\subset
Z'$ is the stratum of $F(z,t)$.
\end{enumerate}
For instance, if $F$ is smooth and equivariant and if there exists a
single stratum $Z'_a$ of $Z'$ such that $F_t(Z\prin)\subset Z'_a$ for
almost all $t$, then $F$ is allowable.

\begin{example}\label{example;radial}
Let $(V,\omega)$ be a symplectic representation space for $G$ as in
Example \ref{example;dilation}.  The radial contraction $F\colon
V\times[\,0,1]\to V$ given by $F(v,t)=tv$ is smooth and equivariant
and satisfies $F_t(Z\prin)\subset Z\prin$ for $\t\ne0$.  Hence it is
allowable.  Likewise, $F$ defines an allowable homotopy
$B\times[\,0,1]\to B$ for any $G$-invariant open ball $B$ about the
origin.
\end{example}

\begin{lemma}\label{lemma;allowable-homotopy}
Let $F\colon M\times[\,0,1]\to M'$ be an allowable homotopy.  Then the
homotopy operator $\kappa_F\colon\Omega(M')\to\Omega(M)$ sends
$\Omega_{\Phi'}(M')$ to $\Omega_\Phi(M)$ and $I_{\Phi'}(M')$ to
$I_\Phi(M)$, and therefore induces a homotopy
$\kappa_F\colon\Omega(X')\to\Omega(X)$.
\end{lemma}

\begin{proof}
Let $\gamma\in\Omega^k_{\Phi'}(M')$.  Then $\gamma$ is invariant, so
$\kappa_F\gamma$ is invariant by \eqref{equation;intertwine}.  Let
$z\in Z\prin$.  Using \eqref{equation;inner} we find that for any
multivector $v\in\Lambda^{k-1}(T_zZ\prin)$
\begin{equation}\label{equation;integral}
i(\xi_M)(\kappa_F\gamma)_z(v)=\int_0^1\phi(t)\,dt,
\end{equation}
where $\phi(t)=-\gamma_{F(z,t)}(\xi_{M'},F_*\partial/\partial
t,(F_t)_*v)$.  Let $Z'_{a(z,t)}$ be the stratum of $Z'$ containing
$F(z,t)$.  Since $F$ is allowable,
$$
F_*\partial/\partial t\in T_{F(z,t)}Z'_{a(z,t)}\qquad\text{and}\qquad
(F_t)_*v\in\Lambda^{k-1}\bigl(T_{F(z,t)}Z'_{a(z,t)}\bigr)
$$
for most $t$.  Moreover, by Lemma
\ref{lemma;restrict}\eqref{item;horizontal} the restriction of
$\gamma$ to $Z'_{a(z,t)}$ is horizontal.  Hence $\phi(t)=0$ for almost
all $t$.  From \eqref{equation;integral} we get
$i(\xi_M)(\kappa_F\gamma)_z(v)=0$; in other words
$\kappa_F\gamma\in\Omega^{k-1}_\Phi(M)$.  The inclusion
$\kappa_FI_{\Phi'}(M')\subset I_\Phi(M)$ is proved in a similar way,
and the last assertion now follows from the isomorphism
\eqref{equation;derham}.
\end{proof}

\begin{example}\label{example;poincare}
Applying Lemma \ref{lemma;allowable-homotopy} to the radial
contraction of Example \ref{example;radial} we find that the de Rham
complex of the symplectic quotient of a vector space $V$ is
homotopically trivial.  More generally, if $Y=(B\cap Z)/G$ is the
symplectic quotient of any $G$-invariant open ball $B$ about the
origin, then the de Rham complex of $Y$ is homotopically trivial.
\end{example}

\begin{example}\label{example;model}
Let $H$ be a closed subgroup of $G$ and let $V$ be a symplectic
$H$-module.  Let $B$ be an $H$-invariant open ball about the origin
and let $O$ be the Hamiltonian $G$-manifold induced by $B$.  Let $Y$
be the symplectic quotient of $B$ by the $H$-action and $U$ the
symplectic quotient of $O$ by the $G$-action.  Then
$\Omega(U)\cong\Omega(Y)$ by Proposition \ref{proposition;induction},
so $\Omega(U)$ is homotopically trivial by Example
\ref{example;poincare}.
\end{example}

This example generalizes to arbitrary Hamiltonian $G$-manifolds by
means of a slice argument.  Let $z\in Z$ and let $H=G_z$ be the
stabilizer of $z$.  Consider the symplectic $H$-module
$V=(T_zGz)^\omega/T_zGz$ known as the \emph{symplectic slice} at $z$.
Choose an $H$-invariant open ball $B$ in $V$ and let $O$ be the
$G$-space induced by $B$.  The symplectic slice theorem due to Marle
and to Guillemin and Sternberg (see e.g.\
\cite[\S~2]{sjamaar-lerman;stratified}) says that, for sufficiently
small $B$, $z$ has a $G$-invariant open neighbourhood that is
isomorphic to $O$ as a Hamiltonian $G$-manifold.  Hence the point
$x\in X$ determined by $z$ has an open neighbourhood $U$ for which
$\Omega(U)$ is homotopically trivial.  By letting $B$ shrink to a
point we obtain a collection of such neighbourhoods, which is a basis
of the topology at $x$.  This finishes the proof of
\eqref{equation;poincare}.

%%%%%%%%%%%%%%%%%%%%%%%%%%%%%%%%%%%%%%%%%%%%%%%%%%%%%%%%%%%%%%%%%%%%%%%%
\section{Integration and the symplectic class}\label{section;class}
%%%%%%%%%%%%%%%%%%%%%%%%%%%%%%%%%%%%%%%%%%%%%%%%%%%%%%%%%%%%%%%%%%%%%%%%

In this section we show that top-degree forms on a compact symplectic
quotient are always integrable and establish a version of Stokes'
theorem.  We conclude that the cohomology class of the symplectic form
and its powers are nonzero.

For technical reasons we do not assume at the outset that $X$ is
compact.  We start by introducing a metric on $X\prin$ and
demonstrating that $X$ has ``locally finite'' volume.  Choose a
$G$-invariant compatible almost complex structure $J$ on the
Hamiltonian $G$-manifold $M$.  The volume element determined by the
Riemannian metric $\sigma=\omega({\cdot},J{\cdot})$ is identical to
the Liouville volume form $\omega^d/d!$ (where $2d=\dim M$).  The
almost complex structure and Riemannian metric descend in a natural
way to each stratum of $X$.  Let $2n=\dim X$ and write
$\mu=\omega\prin^n/n!$ for the volume element of the principal stratum
$X\prin$.

\begin{lemma}\label{lemma;volume}
Every $x\in X$ has an open neighbourhood $U$ such that $\vol U\prin$
is finite.  Hence $X\prin$ has finite volume if $X$ is compact.
\end{lemma}

\begin{proof}
Choose $z\in Z$ mapping to $x$ and let $H=G_z$.  By the symplectic
slice theorem we may take $U$ to be the symplectic quotient of an
$H$-invariant neighbourhood $B$ of the origin in the symplectic slice
$V$ at $z$.  The almost complex structure on $M$ induces one on $V$,
turning $V$ into a unitary $H$-module.  The metric on $U\prin$ induced
by the flat metric $\sigma_V$ on $V$ is quasi-isometric to the metric
induced by $\sigma$.  Therefore it is enough to show that $U$ has
finite volume with respect to the former.  Let $W$ be the orthogonal
complement in $V$ of the subspace of invariants $V^H$.  The quadratic
moment map $\Phi_V$ is constant along $V^H$, so $Z_V=V^H\times Z_W$,
where $Z_V=\Phi_V\inv(0)$ and $Z_W=\Phi_V\inv(0)\cap W$.  Let
$B=B_1\times B_2$, where $B_1$ is an open ball about the origin in
$V^H$ and $B_2$ an $H$-invariant open ball about the origin in $W$.
Then $B$ has a product metric and so do
$(Z_V)\prin=V^H\times(Z_W)\prin$ and the quotient
\begin{equation}\label{equation;product}
U\prin=B_1\times(B_2\cap(Z_W)\prin)/H.
\end{equation}
Recall that the \emph{metric cone} over a Riemannian manifold
$(Y,\sigma_Y)$ is the product $Y\times(0,1)$ with metric
$t^2\sigma_Y+dt\otimes dt$, where $t$ is the coordinate on $(0,1)$.
The metric cone over $Y$ has finite volume if $Y$ does.  For instance,
the ball $B_2$ in $W$ is the metric cone over the sphere $S=\partial
B_2$.  Similarly, with respect to the metric induced by $\sigma_W$,
$B_2\cap(Z_W)\prin$ is a metric cone over $S\cap(Z_W)\prin$.  Upon
taking quotients we see that $(B_2\cap(Z_W)\prin)/H$ is a metric cone
over $(S\cap(Z_W)\prin)/H$.  The link $S\cap Z_W$ is the zero fibre of
the moment map $v\mapsto\bigl(\Phi_W(v),\frac1{2}(1-\abs{v}^2)\bigr)$
for the $H\times\U(1)$-action on $W$, where $\U(1)$ acts by complex
scalar multiplication.  By induction on the depth of the
stratification, the principal stratum of the symplectic quotient
$(S\cap(Z_W))/H$ has finite volume.  Hence $(B_2\cap(Z_W)\prin)/H$ has
finite volume and therefore, because of the product decomposition
\eqref{equation;product}, so does $U\prin$.
\end{proof}

The Riemannian metric on $M$ determines metrics on $\Lambda^k(TM)$ for
all $k$.  Let $\abs{\beta}\in C^0(M)$ denote the pointwise norm of a
form $\beta$ on $M$.  Similarly, for $\alpha\in\Omega(X)$ let
$\abs{\alpha}\in C^0(X\prin)$ denote the pointwise norm over the
principal stratum.  If $\alpha$ is induced by
$\tilde{\alpha}\in\Omega_\Phi(M)$, then $\abs{\tilde{\alpha}}$ is a
$G$-invariant continuous function on $M$ and
\begin{equation}\label{equation;less}
\pi\prin^*\abs{\alpha}\le\iota\prin^*\abs{\tilde{\alpha}}.
\end{equation}
The \emph{support} of a form $\alpha\in\Omega(X)$ is its support as a
section of the sheaf $\Omega$.  This is the same as the closure in $X$
of the support of $\alpha$ considered as a form on $X\prin$.  The
estimate \eqref{equation;less} implies that for $\alpha\in\Omega(X)$
with compact support the pointwise norm $\abs{\alpha}$ is a bounded
function on $X\prin$ and therefore by Lemma \ref{lemma;volume} the
global norm $\int_{X\prin}\abs{\alpha}\mu$ is finite.  In particular,
for $\alpha$ of top degree $2n$ the integral $\int_{X\prin}\alpha$ is
absolutely convergent.

We can now prove Stokes' theorem.  The proof is based on the fact
that the singular strata of $X$ have codimension $\ge2$, which makes
the boundary terms in the integral vanish.

\begin{proposition}\label{proposition;stokes}
$\int_{X\prin}d\gamma=0$ if $\gamma\in\Omega^{2n-1}(X)$ has compact
support.
\end{proposition}

\begin{proof}
We use the notation of the proof of Lemma \ref{lemma;volume}.  By
using partitions of unity we can reduce the general case to the case
where $\gamma$ has compact support in an open subset $U$ of the form
$B_1\times(B_2\cap Z_W)/H$.  Let $2m=\dim Z_W/H$.  If $m=0$ then $U$
is nonsingular and the result follows from the usual version of
Stokes' theorem, so we may assume $m\ge1$.  Let
$\chi\colon[\,0,\infty)\to[\,0,1]$ be a smooth function satisfying
$\chi(t)=0$ for $t$ near $0$ and $\chi(t)=1$ for $t\ge1$.  Define a
sequence of $H$-invariant functions $\tilde{\chi}_k\colon V\to[0,1]$
for $k\ge1$ by $\tilde{\chi}_k(v)=\chi(k\abs{\pr_Wv})$, where
$\pr_W\colon V\to W$ is the orthogonal projection.  These functions
descend to smooth functions $\chi_k\colon U\to[0,1]$.  The functions
$1-\chi_k$ are bump functions supported near the singularities of $U$.
In fact the sets $S_k=\supp(1-\chi_k)$ form a decreasing sequence
satisfying
\begin{equation}\label{equation;singular}
\bigcap_kS_k=B_1\times\{0\mod H\},
\end{equation}
the most singular stratum of $U$.  Therefore $\bigcap_k(S_k)\prin$ is
empty and
$$
\biggl|\int_{X\prin}d\gamma-\int_{X\prin}\chi_kd\gamma\biggr|
=\biggl|\int_{(S_k)\prin}(1-\chi_k)d\gamma\biggr|\le
C\vol(S_k)\prin\to0
$$
as $k\to\infty$.  (Here $C$ is an upper bound for
$\abs{(1-\chi_k)d\gamma}$.)  This shows that
$$
\int_{X\prin}d\gamma=\lim_{k\to\infty}\int_{X\prin}\chi_kd\gamma.
$$
To see that this limit is $0$ we use
$$
\int_{X\prin}\chi_kd\gamma
=\int_{X\prin}d(\chi_k\gamma)-\int_{X\prin}d\chi_k\wedge\gamma.
$$
Since $d(\chi_k\gamma)$ is supported away from the most singular
stratum \eqref{equation;singular}, we can assume by induction on the
depth of the stratification that $\int_{X\prin}d(\chi_k\gamma)=0$.
Moreover,
$$
\biggl|\int_{X\prin}d\chi_k\wedge\gamma\biggr|
\le\int_{X\prin}\abs{d\chi_k}\abs{\gamma}\mu\le
C\int_{(S_k)\prin}\abs{d\chi_k}\mu,
$$
where $C$ is an upper bound for $\abs{\gamma}$.  Let
$\tilde{\rho}_k\colon W\to W$ be the dilation $v\mapsto kv$ and
$\rho_k$ the induced map on $Z_V/H$.  Then $\chi_k=\chi_1\circ\rho_k$
and $S_k=\rho_k\inv(S_1)$.  It follows that
$d\chi_k(x)=kd\chi_1(\rho_k(x))$.  By \eqref{equation;product},
$U\prin$ is the product of a ball and a metric cone, so
$\vol(S_k)\prin=k^{-2m}\vol(S_1)\prin$, where $2m=\dim Z_W/H\ge2$.
Hence
$$
\biggl|\int_{X\prin}d\chi_k\wedge\gamma\biggr|\le
Ck^{1-2m}\int_{(S_1)\prin}\abs{d\chi_1}\mu\to0
$$
as $k\to\infty$.  Therefore
$\lim_{k\to\infty}\int_{X\prin}\chi_kd\gamma=0$.
\end{proof}

Stokes' theorem implies that the volume form of a compact quotient is
not exact.

\begin{corollary}
Suppose that $X$ is compact.  Then the class of $\omega\prin^k$ in
$H^{2k}(\Omega(X))$ is nonzero for $0\le k\le n$, where $2n=\dim X$.
\end{corollary}

%%%%%%%%%%%%%%%%%%%%%%%%%%%%%%%%%%%%%%%%%%%%%%%%%%%%%%%%%%%%%%%%%%%%%%%%
\section{Addendum}\label{section;addendum}
%%%%%%%%%%%%%%%%%%%%%%%%%%%%%%%%%%%%%%%%%%%%%%%%%%%%%%%%%%%%%%%%%%%%%%%%

The above results can be generalized in two obvious ways.  First we
consider symplectic quotients at nonzero levels.  Let $\ca{O}$ be a
coadjoint orbit in $\g^*$.  The \emph{symplectic quotient} at $\ca{O}$
is $X_\ca{O}=Z_\ca{O}/G$, where $Z_\ca{O}$ is the fibre
$\Phi\inv(\ca{O})$.  The spaces $Z_\ca{O}$ and $X_\ca{O}$ stratify in
exactly the same way as when $\ca{O}=\{0\}$ and the strata of
$X_\ca{O}$ again carry natural symplectic forms.  Differential forms
on $X_\ca{O}$ can now be defined as before.  There is a symplectic
slice theorem for orbits in $Z_\ca{O}$, so all our results generalize
to this situation with virtually unchanged proofs.

Next we consider actions of a noncompact group $G$.  The symplectic
slice theorem remains valid, provided that $G$ acts properly on $M$.
For locally closed coadjoint orbits $\ca{O}$ stratifications of
$Z_\ca{O}$ and $X_\ca{O}$ were obtained in \cite{bates-lerman;proper}.
However, our definition of forms on $X$ is valid as it stands only
when $\ca{O}$ is closed, because forms on a nonclosed subset may not
extend to the ambient manifold.  If $\ca{O}$ is locally closed we
define $\Omega(X_\ca{O})=\Omega_\Phi(N)/I_\Phi(N)$.  Here
$N=\Phi\inv(D)$ is the preimage of any $G$-invariant open
neighbourhood $D$ of $\ca{O}$ in $\g^*$ such that $\ca{O}$ is closed
in $D$, $\Omega_\Phi(N)$ is the set of $G$-invariant forms on $N$ that
restrict to basic forms on $(Z_\ca{O})\prin$, and $I_\Phi(N)$ is the
set of $G$-invariant forms on $N$ that restrict to $0$ on
$(Z_\ca{O})\prin$.  With this minor modification our results carry
over to symplectic quotients by proper actions at locally closed
coadjoint orbits.  (For general orbits one might try to apply the
methods developed in \cite{cushman-sniatycki;differential-structure},
but we have not attempted this.)

%%%%%%%%%%%%%%%%%%%%%%%%%%%%%%%%%%%%%%%%%%%%%%%%%%%%%%%%%%%%%%%%%%%%%%%%
\appendix
%%%%%%%%%%%%%%%%%%%%%%%%%%%%%%%%%%%%%%%%%%%%%%%%%%%%%%%%%%%%%%%%%%%%%%%%

%%%%%%%%%%%%%%%%%%%%%%%%%%%%%%%%%%%%%%%%%%%%%%%%%%%%%%%%%%%%%%%%%%%%%%%%
\section{Forms on homogeneous bundles}\label{section;homogeneous}
%%%%%%%%%%%%%%%%%%%%%%%%%%%%%%%%%%%%%%%%%%%%%%%%%%%%%%%%%%%%%%%%%%%%%%%%

Let $G$ be a compact Lie group and $H$ a closed subgroup.  For any
$H$-manifold $F$ we can form the homogeneous fibre bundle with fibre
$F$ over $G/H$,
$$
E=(G\times F)/H.
$$
The map $f\colon F\to E$ defined by $f(p)=[1,p]$ identifies $F$ with
the fibre over the coset $0\mod H$.  (Here $[g,p]$ denotes the coset
$(g,p)\mod H$ of $(g,p)\in G\times F$.)  Restriction to the fibre is a
homomorphism
$$
f^*\colon\Omega(E)^G\to\Omega(F)^H.
$$
It is not hard to see that $G$-basic forms on $E$ restrict to
$H$-basic forms on $F$ and that
$f^*\colon\Omega\bas(E)\to\Omega\bas(F)$ is an isomorphism.  We
require a slight generalization of this elementary fact.

Choose an $H$-equivariant projection $\g\to\h$; this determines a
$G$-invariant connection $1$-form $\theta\in\Omega^1(G,\h)^{G\times
H}$ on the principal $H$-bundle $G\to G/H$.  Let $VE$ be the vertical
tangent bundle of $E\to G/H$ and let $\theta_E\in\Omega^1(E,VE)^G$ be
the $G$-invariant connection $1$-form on $E$ associated to $\theta$.
Let $\gamma\in\Omega(F)^H$ be any invariant form on the fibre.  Define
a form $e(\gamma)\in\Omega(E)$ by putting
$$
e(\gamma)_{[g,p]}(v)=\gamma_p\bigl((g\inv)_*\theta_E(v)\bigr)
$$
for $[g,p]\in E$ and $v\in\Lambda(T_{[g,p]}E)$.  (For simplicity we
write $\theta_E$ for the extension of the connection $\theta_E\colon
TE\to VE$ to a multiplicative map $\Lambda(TE)\to\Lambda(VE)$.)  The
$H$-invariance of $\gamma$ implies that $e(\gamma)_{[g,p]}(v)$ does
not depend on the choice of the representative $(g,p)$ of the coset
$[g,p]$.  The $G$-invariance of $\theta_E$ implies that $e(\gamma)$ is
$G$-invariant.  Thus we have defined a map
$$
e\colon\Omega(F)^H\to\Omega(E)^G,
$$
which we call the \emph{extension homomorphism} determined by
$\theta$.  (An alternative definition runs as follows.  Let
$\ca{V}=\pr^*TF$, where $\pr\colon G\times F\to F$ is the Cartesian
projection.  The vertical bundle of $E$ is then the quotient
$VE\cong\ca{V}/H$.  A form $\gamma\in\Omega(F)^H$ is a section of
$\Lambda(TF)$ and as such extends uniquely to a section
$\tilde{\gamma}$ of $\ca{V}$ that is constant along $G$.  Then
$\tilde{\gamma}$ is $G\times H$-invariant and so descends to a
$G$-invariant section $\bar{\gamma}$ of $VE$.  Thus
$e(\gamma)=\theta_E^*\bar{\gamma}$ is a $G$-invariant section of
$\Lambda(TE)$.  This argument also shows that $e(\gamma)$ is smooth.)
The following result is immediate from the definition.

\begin{lemma}\label{lemma;extension}
\begin{enumerate}
\item\label{item;right}
$f^*e(\gamma)=\gamma$ for $\gamma\in\Omega(F)^H$;
\item\label{item;basic}
$e$ maps $\Omega(F)\bas$ to $\Omega(E)\bas$;
\item\label{item;left}
$e(f^*\beta)=\beta$ for $\beta\in\Omega(E)\bas$.
\end{enumerate}
\end{lemma}

It follows from \eqref{item;right} that
$f^*\colon\Omega(E)^G\to\Omega(F)^H$ is surjective and from
\eqref{item;basic}--\eqref{item;left} that
$f^*\colon\Omega\bas(E)\to\Omega\bas(F)$ is an isomorphism, as noted
above.  Now let $F'$ be a second $H$-manifold and let $j\colon F\to
F'$ be an $H$-equivariant map.  Then $j$ extends naturally to an
$G$-equivariant bundle map $\bar{\jmath}\colon E\to E'=(G\times
F')/H$.  Moreover $\theta_E$ is the pullback of the associated
connection $\theta_{E'}$ on $E'$.  This implies that the extension
homomorphism is functorial in the following sense.

\begin{lemma}\label{lemma;functor}
$e\circ j^*=\bar{\jmath}^*\circ e'$, where
$e'\colon\Omega(F')^H\to\Omega(E')^G$ is the extension homomorphism
for $E'$.
\end{lemma}

%%%%%%%%%%%%%%%%%%%%%%%%%%%%%%%%%%%%%%%%%%%%%%%%%%%%%%%%%%%%%%%%%%%%%%%%
%%%%%%%%%%%%%%%%%%%%%%%%%%%%%%%%%%%%%%%%%%%%%%%%%%%%%%%%%%%%%%%%%%%%%%%%

\providecommand{\bysame}{\leavevmode\hbox to3em{\hrulefill}\thinspace}
\providecommand{\MR}{\relax\ifhmode\unskip\space\fi MR }
\providecommand{\MRhref}[2]{%
  \href{http://www.ams.org/mathscinet-getitem?mr=#1}{#2}
}
\providecommand{\href}[2]{#2}

%%%%%%%%%%%%%%%%%%%%%%%%%%%%%%%%%%%%%%%%%%%%%%%%%%%%%%%%%%%%%%%%%%%%%%%%
%%%%%%%%%%%%%%%%%%%%%%%%%%%%%%%%%%%%%%%%%%%%%%%%%%%%%%%%%%%%%%%%%%%%%%%%

\end{document}